\catcode`\@=11
\font\tensmc=cmcsc10      
\def\smc{\tensmc}

\def\hcorrection#1{\advance\hoffset by #1 }
\def\vcorrection#1{\advance\voffset by #1 }
\def\wlog#1{}
\newif\iftitle@
\outer\def\title{\title@true\vglue 24\p@ plus 12\p@ minus 12\p@
   \bgroup\let\\=\cr\tabskip\centering
   \halign to \hsize\bgroup\tenbf\hfill\ignorespaces##\unskip\hfill\cr}
\def\endtitle{\cr\egroup\egroup\vglue 18\p@ plus 12\p@ minus 6\p@}
\outer\def\author{\iftitle@\vglue -18\p@ plus -12\p@ minus
-6\p@\fi\vglue
    12\p@ plus 6\p@ minus 3\p@\bgroup\let\\=\cr\tabskip\centering
    \halign to \hsize\bgroup\smc\hfill\ignorespaces##\unskip\hfill\cr}
\def\endauthor{\cr\egroup\egroup\vglue 18\p@ plus 12\p@ minus 6\p@}
\outer\def\heading{\bigbreak\bgroup\let\\=\cr\tabskip\centering
    \halign to \hsize\bgroup\smc\hfill\ignorespaces##\unskip\hfill\cr}
\def\endheading{\cr\egroup\egroup\nobreak\medskip}
\outer\def\subheading#1{\medbreak\noindent{\tenbf\ignorespaces
      #1\unskip.\enspace}\ignorespaces}
\outer\def\proclaim#1{\medbreak\noindent\smc\ignorespaces
    #1\unskip.\enspace\sl\ignorespaces}
\outer\def\endproclaim{\par\ifdim\lastskip<\medskipamount\removelastskip
  \penalty 55 \fi\medskip\rm}
\outer\def\demo#1{\par\ifdim\lastskip<\smallskipamount\removelastskip
    \smallskip\fi\noindent{\smc\ignorespaces#1\unskip:\enspace}\rm
      \ignorespaces}

\newcount\footmarkcount@
\footmarkcount@=1
\def\makefootnote@#1#2{\insert\footins{\interlinepenalty=100
  \splittopskip=\ht\strutbox \splitmaxdepth=\dp\strutbox
  \floatingpenalty=\@MM
  \leftskip=\z@\rightskip=\z@\spaceskip=\z@\xspaceskip=\z@
  \noindent{#1}\footstrut\rm\ignorespaces #2\strut}}
\def\footnote{\let\@sf=\empty\ifhmode\edef\@sf{\spacefactor
   =\the\spacefactor}\/\fi\futurelet\next\footnote@}
\def\footnote@{\ifx"\next\let\next\footnote@@\else
    \let\next\footnote@@@\fi\next}
\def\footnote@@"#1"#2{#1\@sf\relax\makefootnote@{#1}{#2}}
\def\footnote@@@#1{$^{\number\footmarkcount@}$\makefootnote@
   {$^{\number\footmarkcount@}$}{#1}\global\advance\footmarkcount@ by 1 }

\hyphenation{man-u-script man-u-scripts ap-pen-dix ap-pen-di-ces}
\hyphenation{data-base data-bases}
\ifx\amstexloaded@\relax\catcode`\@=13
  \endinput\else\let\amstexloaded@=\relax\fi
\newlinechar=`\^^J
\def\eat@#1{}
\def\Space@.{\futurelet\Space@\relax}
\Space@. %
\newhelp\athelp@
{Only certain combinations beginning with @ make sense to me.^^J
Perhaps you wanted \string\@\space for a printed @?^^J I've
ignored the character or group after @.}
\def\futureletnextat@{\futurelet\next\at@}
{\catcode`\@=\active \lccode`\Z=`\@ \lowercase
{\gdef@{\expandafter\csname futureletnextatZ\endcsname}
\expandafter\gdef\csname atZ\endcsname
   {\ifcat\noexpand\next a\def\next{\csname atZZ\endcsname}\else
   \ifcat\noexpand\next0\def\next{\csname atZZ\endcsname}\else
    \def\next{\csname atZZZ\endcsname}\fi\fi\next}
\expandafter\gdef\csname atZZ\endcsname#1{\expandafter
   \ifx\csname #1Zat\endcsname\relax\def\next
     {\errhelp\expandafter=\csname athelpZ\endcsname
      \errmessage{Invalid use of \string@}}\else
       \def\next{\csname #1Zat\endcsname}\fi\next}
\expandafter\gdef\csname atZZZ\endcsname#1{\errhelp
    \expandafter=\csname athelpZ\endcsname
      \errmessage{Invalid use of \string@}}}}
\def\atdef@#1{\expandafter\def\csname #1@at\endcsname}
\newhelp\defahelp@{If you typed \string\define\space cs instead of
\string\define\string\cs\space^^J I've substituted an inaccessible
control sequence so that your^^J definition will be completed
without mixing me up too badly.^^J If you typed
\string\define{\string\cs} the inaccessible control sequence^^J
was defined to be \string\cs, and the rest of your^^J definition
appears as input.}
\newhelp\defbhelp@{I've ignored your definition, because it might^^J
conflict with other uses that are important to me.}
\def\define{\futurelet\next\define@}
\def\define@{\ifcat\noexpand\next\relax
  \def\next{\define@@}%
  \else\errhelp=\defahelp@
  \errmessage{\string\define\space must be followed by a control
     sequence}\def\next{\def\garbage@}\fi\next}
\def\undefined@{}
\def\preloaded@{}
\def\define@@#1{\ifx#1\relax\errhelp=\defbhelp@
   \errmessage{\string#1\space is already defined}\def\next{\def\garbage@}%
   \else\expandafter\ifx\csname\expandafter\eat@\string
     #1@\endcsname\undefined@\errhelp=\defbhelp@
   \errmessage{\string#1\space can't be defined}\def\next{\def\garbage@}%
   \else\expandafter\ifx\csname\expandafter\eat@\string#1\endcsname\relax
     \def\next{\def#1}\else\errhelp=\defbhelp@
     \errmessage{\string#1\space is already defined}\def\next{\def\garbage@}%
      \fi\fi\fi\next}
\def\famzero{\fam\z@}

\def\inf{\mathop{\famzero inf}}

\def\lim{\mathop{\famzero lim}}

\def\max{\mathop{\famzero max}}
\def\min{\mathop{\famzero min}}

\def\textfont@#1#2{\def#1{\relax\ifmmode
    \errmessage{Use \string#1\space only in text}\else#2\fi}}
\textfont@\rm\tenrm \textfont@\it\tenit \textfont@\sl\tensl
\textfont@\bf\tenbf \textfont@\smc\tensmc
\let\ic@=\/
\def\/{\unskip\ic@}
\def\textfonti{\the\textfont1 }
\def\t#1#2{{\edef\next{\the\font}\textfonti\accent"7F \next#1#2}}
\let\B=\=
\let\D=\.
\def~{\unskip\nobreak\ \ignorespaces}
{\catcode`\@=\active \gdef\@{\char'100 }}
\atdef@-{\leavevmode\futurelet\next\athyph@}
\def\athyph@{\ifx\next-\let\next=\athyph@@
  \else\let\next=\athyph@@@\fi\next}
\def\athyph@@@{\hbox{-}}
\def\athyph@@#1{\futurelet\next\athyph@@@@}
\def\athyph@@@@{\if\next-\def\next##1{\hbox{---}}\else
    \def\next{\hbox{--}}\fi\next}
\def\.{.\spacefactor=\@m}
\atdef@.{\null.} \atdef@,{\null,} \atdef@;{\null;}
\atdef@:{\null:} \atdef@?{\null?} \atdef@!{\null!}
\def\srdr@{\thinspace}
\def\drsr@{\kern.02778em}
\def\sldl@{\kern.02778em}
\def\dlsl@{\thinspace}
\atdef@"{\unskip\futurelet\next\atqq@}
\def\atqq@{\ifx\next\Space@\def\next. {\atqq@@}\else
     \def\next.{\atqq@@}\fi\next.}
\def\atqq@@{\futurelet\next\atqq@@@}
\def\atqq@@@{\ifx\next`\def\next`{\atqql@}\else\def\next'{\atqqr@}\fi\next}
\def\atqql@{\futurelet\next\atqql@@}
\def\atqql@@{\ifx\next`\def\next`{\sldl@``}\else\def\next{\dlsl@`}\fi\next}
\def\atqqr@{\futurelet\next\atqqr@@}
\def\atqqr@@{\ifx\next'\def\next'{\srdr@''}\else\def\next{\drsr@'}\fi\next}
\def\flushpar{\par\noindent}
\def\textfontii{\the\textfont2 }
\def\{{\relax\ifmmode\lbrace\else
    {\textfontii f}\spacefactor=\@m\fi}
\def\}{\relax\ifmmode\rbrace\else
    \let\@sf=\empty\ifhmode\edef\@sf{\spacefactor=\the\spacefactor}\fi
      {\textfontii g}\@sf\relax\fi}
\def\nonhmodeerr@#1{\errmessage
     {\string#1\space allowed only within text}}
\def\linebreak{\relax\ifhmode\unskip\break\else
    \nonhmodeerr@\linebreak\fi}
\def\allowlinebreak{\relax
   \ifhmode\allowbreak\else\nonhmodeerr@\allowlinebreak\fi}
\newskip\saveskip@
\def\nolinebreak{\relax\ifhmode\saveskip@=\lastskip\unskip
  \nobreak\ifdim\saveskip@>\z@\hskip\saveskip@\fi
   \else\nonhmodeerr@\nolinebreak\fi}
\def\newline{\relax\ifhmode\null\hfil\break
    \else\nonhmodeerr@\newline\fi}
\def\nonmathaerr@#1{\errmessage
     {\string#1\space is not allowed in display math mode}}
\def\nonmathberr@#1{\errmessage{\string#1\space is allowed only in math mode}}
\def\mathbreak{\relax\ifmmode\ifinner\break\else
   \nonmathaerr@\mathbreak\fi\else\nonmathberr@\mathbreak\fi}
\def\nomathbreak{\relax\ifmmode\ifinner\nobreak\else
    \nonmathaerr@\nomathbreak\fi\else\nonmathberr@\nomathbreak\fi}
\def\allowmathbreak{\relax\ifmmode\ifinner\allowbreak\else
     \nonmathaerr@\allowmathbreak\fi\else\nonmathberr@\allowmathbreak\fi}
\def\pagebreak{\relax\ifmmode
   \ifinner\errmessage{\string\pagebreak\space
     not allowed in non-display math mode}\else\postdisplaypenalty-\@M\fi
   \else\ifvmode\penalty-\@M\else\edef\spacefactor@
       {\spacefactor=\the\spacefactor}\vadjust{\penalty-\@M}\spacefactor@
    \relax\fi\fi}
\def\nopagebreak{\relax\ifmmode
     \ifinner\errmessage{\string\nopagebreak\space
    not allowed in non-display math mode}\else\postdisplaypenalty\@M\fi
    \else\ifvmode\nobreak\else\edef\spacefactor@
    {\spacefactor=\the\spacefactor}\vadjust{\penalty\@M}\spacefactor@
     \relax\fi\fi}
\def\newpage{\relax\ifvmode\vfill\penalty-\@M\else\nonvmodeerr@\newpage\fi}
\def\nonvmodeerr@#1{\errmessage
    {\string#1\space is allowed only between paragraphs}}
\def\smallpagebreak{\relax\ifvmode\smallbreak
      \else\nonvmodeerr@\smallpagebreak\fi}
\def\medpagebreak{\relax\ifvmode\medbreak
       \else\nonvmodeerr@\medpagebreak\fi}
\def\bigpagebreak{\relax\ifvmode\bigbreak
      \else\nonvmodeerr@\bigpagebreak\fi}
\newdimen\captionwidth@
\captionwidth@=\hsize \advance\captionwidth@ by -1.5in
\def\caption#1{}
\def\topspace#1{\gdef\thespace@{#1}\ifvmode\def\next
    {\futurelet\next\topspace@}\else\def\next{\nonvmodeerr@\topspace}\fi\next}
\def\topspace@{\ifx\next\Space@\def\next. {\futurelet\next\topspace@@}\else
     \def\next.{\futurelet\next\topspace@@}\fi\next.}
\def\topspace@@{\ifx\next\caption\let\next\topspace@@@\else
    \let\next\topspace@@@@\fi\next}
 \def\topspace@@@@{\topinsert\vbox to
       \thespace@{}\endinsert}
\def\topspace@@@\caption#1{\topinsert\vbox to
    \thespace@{}\nobreak
      \smallskip
    \setbox\z@=\hbox{\noindent\ignorespaces#1\unskip}%
   \ifdim\wd\z@>\captionwidth@
   \centerline{\vbox{\hsize=\captionwidth@\noindent\ignorespaces#1\unskip}}%
   \else\centerline{\box\z@}\fi\endinsert}
\def\midspace#1{\gdef\thespace@{#1}\ifvmode\def\next
    {\futurelet\next\midspace@}\else\def\next{\nonvmodeerr@\midspace}\fi\next}
\def\midspace@{\ifx\next\Space@\def\next. {\futurelet\next\midspace@@}\else
     \def\next.{\futurelet\next\midspace@@}\fi\next.}
\def\midspace@@{\ifx\next\caption\let\next\midspace@@@\else
    \let\next\midspace@@@@\fi\next}
 \def\midspace@@@@{\midinsert\vbox to
       \thespace@{}\endinsert}
\def\midspace@@@\caption#1{\midinsert\vbox to
    \thespace@{}\nobreak
      \smallskip
      \setbox\z@=\hbox{\noindent\ignorespaces#1\unskip}%
      \ifdim\wd\z@>\captionwidth@
    \centerline{\vbox{\hsize=\captionwidth@\noindent\ignorespaces#1\unskip}}%
    \else\centerline{\box\z@}\fi\endinsert}
\mathchardef\prime@="0230
\def\prime{{{}\prime@{}}}
\def\prim@s{\prime@\futurelet\next\pr@m@s}

\def\,{\relax\ifmmode\mskip\thinmuskip\else\thinspace\fi}
\def\!{\relax\ifmmode\mskip-\thinmuskip\else\negthinspace\fi}
\def\frac#1#2{{#1\over#2}}

\def\:{\nobreak\hskip.1111em{:}\hskip.3333em plus .0555em\relax}
\def\intic@{\mathchoice{\hskip5\p@}{\hskip4\p@}{\hskip4\p@}{\hskip4\p@}}
\def\negintic@
 {\mathchoice{\hskip-5\p@}{\hskip-4\p@}{\hskip-4\p@}{\hskip-4\p@}}
\def\intkern@{\mathchoice{\!\!\!}{\!\!}{\!\!}{\!\!}}
\def\intdots@{\mathchoice{\cdots}{{\cdotp}\mkern1.5mu
    {\cdotp}\mkern1.5mu{\cdotp}}{{\cdotp}\mkern1mu{\cdotp}\mkern1mu
      {\cdotp}}{{\cdotp}\mkern1mu{\cdotp}\mkern1mu{\cdotp}}}
\newcount\intno@
\def\iint{\intno@=\tw@\futurelet\next\ints@}
\def\iiint{\intno@=\thr@@\futurelet\next\ints@}
\def\iiiint{\intno@=4 \futurelet\next\ints@}
\def\idotsint{\intno@=\z@\futurelet\next\ints@}
\def\ints@{\findlimits@\ints@@}
\newif\iflimtoken@
\newif\iflimits@
\def\findlimits@{\limtoken@false\limits@false\ifx\next\limits
 \limtoken@true\limits@true\else\ifx\next\nolimits\limtoken@true\limits@false
    \fi\fi}
\def\multintlimits@{\intop\ifnum\intno@=\z@\intdots@
  \else\intkern@\fi
    \ifnum\intno@>\tw@\intop\intkern@\fi
     \ifnum\intno@>\thr@@\intop\intkern@\fi\intop}
\def\multint@{\int\ifnum\intno@=\z@\intdots@\else\intkern@\fi
   \ifnum\intno@>\tw@\int\intkern@\fi
    \ifnum\intno@>\thr@@\int\intkern@\fi\int}
\def\ints@@{\iflimtoken@\def\ints@@@{\iflimits@
   \negintic@\mathop{\intic@\multintlimits@}\limits\else
    \multint@\nolimits\fi\eat@}\else
     \def\ints@@@{\multint@\nolimits}\fi\ints@@@}
\def\Sb{_\bgroup\vspace@
    \baselineskip=\fontdimen10 \scriptfont\tw@
    \advance\baselineskip by \fontdimen12 \scriptfont\tw@
    \lineskip=\thr@@\fontdimen8 \scriptfont\thr@@
    \lineskiplimit=\thr@@\fontdimen8 \scriptfont\thr@@
    \Let@\vbox\bgroup\halign\bgroup \hfil$\scriptstyle
        {##}$\hfil\cr}
\def\endSb{\crcr\egroup\egroup\egroup}
\def\Sp{^\bgroup\vspace@
    \baselineskip=\fontdimen10 \scriptfont\tw@
    \advance\baselineskip by \fontdimen12 \scriptfont\tw@
    \lineskip=\thr@@\fontdimen8 \scriptfont\thr@@
    \lineskiplimit=\thr@@\fontdimen8 \scriptfont\thr@@
    \Let@\vbox\bgroup\halign\bgroup \hfil$\scriptstyle
        {##}$\hfil\cr}
\def\endSp{\crcr\egroup\egroup\egroup}
\def\Let@{\relax\iffalse{\fi\let\\=\cr\iffalse}\fi}
\def\vspace@{\def\vspace##1{\noalign{\vskip##1 }}}
\def\aligned{\,\vcenter\bgroup\vspace@\Let@\openup\jot\m@th\ialign
  \bgroup \strut\hfil$\displaystyle{##}$&$\displaystyle{{}##}$\hfil\crcr}
\def\endaligned{\crcr\egroup\egroup}
\def\matrix{\,\vcenter\bgroup\Let@\vspace@
    \normalbaselines
  \m@th\ialign\bgroup\hfil$##$\hfil&&\quad\hfil$##$\hfil\crcr
    \mathstrut\crcr\noalign{\kern-\baselineskip}}
\def\endmatrix{\crcr\mathstrut\crcr\noalign{\kern-\baselineskip}\egroup
        \egroup\,}
\newtoks\hashtoks@
\hashtoks@={#}
\def\format{\crcr\egroup\iffalse{\fi\ifnum`}=0 \fi\format@}
\def\format@#1\\{\def\preamble@{#1}%
  \def\c{\hfil$\the\hashtoks@$\hfil}%
  \def\r{\hfil$\the\hashtoks@$}%
  \def\l{$\the\hashtoks@$\hfil}%
  \setbox\z@=\hbox{\xdef\Preamble@{\preamble@}}\ifnum`{=0 \fi\iffalse}\fi
   \ialign\bgroup\span\Preamble@\crcr}

\def\cases{\left\{\,\vcenter\bgroup\vspace@
     \normalbaselines\openup\jot\m@th
       \Let@\ialign\bgroup$##$\hfil&\quad$##$\hfil\crcr
      \mathstrut\crcr\noalign{\kern-\baselineskip}}
\def\endcases{\endmatrix\right.}
\newif\iftagsleft@
\tagsleft@true
\def\TagsOnRight{\global\tagsleft@false}
\def\tag#1$${\iftagsleft@\leqno\else\eqno\fi
 \hbox{\def\pagebreak{\global\postdisplaypenalty-\@M}%
 \def\nopagebreak{\global\postdisplaypenalty\@M}\rm(#1\unskip)}%
  $$\postdisplaypenalty\z@\ignorespaces}
\interdisplaylinepenalty=\@M
\def\allowdisplaybreak@{\def\allowdisplaybreak{\noalign{\allowbreak}}}
\def\displaybreak@{\def\displaybreak{\noalign{\break}}}
\def\align#1\endalign{\def\tag{&}\vspace@\allowdisplaybreak@\displaybreak@
  \iftagsleft@\lalign@#1\endalign\else
   \ralign@#1\endalign\fi}
\def\ralign@#1\endalign{\displ@y\Let@\tabskip\centering\halign to\displaywidth
     {\hfil$\displaystyle{##}$\tabskip=\z@&$\displaystyle{{}##}$\hfil
       \tabskip=\centering&\llap{\hbox{(\rm##\unskip)}}\tabskip\z@\crcr
         #1\crcr}}
\def\lalign@
 #1\endalign{\displ@y\Let@\tabskip\centering\halign to \displaywidth
   {\hfil$\displaystyle{##}$\tabskip=\z@&$\displaystyle{{}##}$\hfil
   \tabskip=\centering&\kern-\displaywidth
    \rlap{\hbox{(\rm##\unskip)}}\tabskip=\displaywidth\crcr
           #1\crcr}}
\def\overrightarrow{\mathpalette\overrightarrow@}
\def\overrightarrow@#1#2{\vbox{\ialign{$##$\cr
    #1{-}\mkern-6mu\cleaders\hbox{$#1\mkern-2mu{-}\mkern-2mu$}\hfill
     \mkern-6mu{\to}\cr
     \noalign{\kern -1\p@\nointerlineskip}
     \hfil#1#2\hfil\cr}}}
\def\overleftarrow{\mathpalette\overleftarrow@}
\def\overleftarrow@#1#2{\vbox{\ialign{$##$\cr
     #1{\leftarrow}\mkern-6mu\cleaders\hbox{$#1\mkern-2mu{-}\mkern-2mu$}\hfill
      \mkern-6mu{-}\cr
     \noalign{\kern -1\p@\nointerlineskip}
     \hfil#1#2\hfil\cr}}}
\def\overleftrightarrow{\mathpalette\overleftrightarrow@}
\def\overleftrightarrow@#1#2{\vbox{\ialign{$##$\cr
     #1{\leftarrow}\mkern-6mu\cleaders\hbox{$#1\mkern-2mu{-}\mkern-2mu$}\hfill
       \mkern-6mu{\to}\cr
    \noalign{\kern -1\p@\nointerlineskip}
      \hfil#1#2\hfil\cr}}}
\def\underrightarrow{\mathpalette\underrightarrow@}
\def\underrightarrow@#1#2{\vtop{\ialign{$##$\cr
    \hfil#1#2\hfil\cr
     \noalign{\kern -1\p@\nointerlineskip}
    #1{-}\mkern-6mu\cleaders\hbox{$#1\mkern-2mu{-}\mkern-2mu$}\hfill
     \mkern-6mu{\to}\cr}}}
\def\underleftarrow{\mathpalette\underleftarrow@}
\def\underleftarrow@#1#2{\vtop{\ialign{$##$\cr
     \hfil#1#2\hfil\cr
     \noalign{\kern -1\p@\nointerlineskip}
     #1{\leftarrow}\mkern-6mu\cleaders\hbox{$#1\mkern-2mu{-}\mkern-2mu$}\hfill
      \mkern-6mu{-}\cr}}}
\def\underleftrightarrow{\mathpalette\underleftrightarrow@}
\def\underleftrightarrow@#1#2{\vtop{\ialign{$##$\cr
      \hfil#1#2\hfil\cr
    \noalign{\kern -1\p@\nointerlineskip}
     #1{\leftarrow}\mkern-6mu\cleaders\hbox{$#1\mkern-2mu{-}\mkern-2mu$}\hfill
       \mkern-6mu{\to}\cr}}}
\def\sqrt#1{\radical"270370 {#1}}
\def\dots{\relax\ifmmode\let\next=\ldots\else\let\next=\tdots@\fi\next}
\def\tdots@{\unskip\ \tdots@@}
\def\tdots@@{\futurelet\next\tdots@@@}
\def\tdots@@@{$\mathinner{\ldotp\ldotp\ldotp}\,
   \ifx\next,$\else
   \ifx\next.\,$\else
   \ifx\next;\,$\else
   \ifx\next:\,$\else
   \ifx\next?\,$\else
   \ifx\next!\,$\else
   $ \fi\fi\fi\fi\fi\fi}
\def\text{\relax\ifmmode\let\next=\text@\else\let\next=\text@@\fi\next}
\def\text@@#1{\hbox{#1}}
\def\text@#1{\mathchoice
 {\hbox{\everymath{\displaystyle}\def\textfonti{\the\textfont1 }%
    \def\textfontii{\the\textfont2 }\textdef@@ T#1}}
 {\hbox{\everymath{\textstyle}\def\textfonti{\the\textfont1 }%
    \def\textfontii{\the\textfont2 }\textdef@@ T#1}}
 {\hbox{\everymath{\scriptstyle}\def\textfonti{\the\scriptfont1 }%
   \def\textfontii{\the\scriptfont2 }\textdef@@ S\rm#1}}
 {\hbox{\everymath{\scriptscriptstyle}\def\textfonti{\the\scriptscriptfont1 }%
   \def\textfontii{\the\scriptscriptfont2 }\textdef@@ s\rm#1}}}
\def\textdef@@#1{\textdef@#1\rm \textdef@#1\bf
   \textdef@#1\sl \textdef@#1\it}

\def\textdef@#1#2{\def\next{\csname\expandafter\eat@\string#2fam\endcsname}%
\if S#1\edef#2{\the\scriptfont\next\relax}%
 \else\if s#1\edef#2{\the\scriptscriptfont\next\relax}%
 \else\edef#2{\the\textfont\next\relax}\fi\fi}
\scriptfont\itfam=\tenit \scriptscriptfont\itfam=\tenit
\scriptfont\slfam=\tensl \scriptscriptfont\slfam=\tensl
\mathcode`\0="0030 \mathcode`\1="0031 \mathcode`\2="0032
\mathcode`\3="0033 \mathcode`\4="0034 \mathcode`\5="0035
\mathcode`\6="0036 \mathcode`\7="0037 \mathcode`\8="0038
\mathcode`\9="0039
\def\Cal{\relax\ifmmode\let\next=\Cal@\else
     \def\next{\errmessage{Use \string\Cal\space only in math mode}}\fi\next}
\def\Cal@#1{{\fam2 #1}}
\def\bold{\relax\ifmmode\let\next=\bold@\else
   \def\next{\errmessage{Use \string\bold\space only in math
      mode}}\fi\next}\def\bold@#1{{\fam\bffam #1}}
\mathchardef\Gamma="0000 \mathchardef\Delta="0001
\mathchardef\Theta="0002 \mathchardef\Lambda="0003
\mathchardef\Xi="0004 \mathchardef\Pi="0005
\mathchardef\Sigma="0006 \mathchardef\Upsilon="0007
\mathchardef\Phi="0008 \mathchardef\Psi="0009
\mathchardef\Omega="000A \mathchardef\varGamma="0100
\mathchardef\varDelta="0101 \mathchardef\varTheta="0102
\mathchardef\varLambda="0103 \mathchardef\varXi="0104
\mathchardef\varPi="0105 \mathchardef\varSigma="0106
\mathchardef\varUpsilon="0107 \mathchardef\varPhi="0108
\mathchardef\varPsi="0109 \mathchardef\varOmega="010A
\font\dummyft@=dummy \fontdimen1 \dummyft@=\z@ \fontdimen2
\dummyft@=\z@ \fontdimen3 \dummyft@=\z@ \fontdimen4 \dummyft@=\z@
\fontdimen5 \dummyft@=\z@ \fontdimen6 \dummyft@=\z@ \fontdimen7
\dummyft@=\z@ \fontdimen8 \dummyft@=\z@ \fontdimen9 \dummyft@=\z@
\fontdimen10 \dummyft@=\z@ \fontdimen11 \dummyft@=\z@ \fontdimen12
\dummyft@=\z@ \fontdimen13 \dummyft@=\z@ \fontdimen14
\dummyft@=\z@ \fontdimen15 \dummyft@=\z@ \fontdimen16
\dummyft@=\z@ \fontdimen17 \dummyft@=\z@ \fontdimen18
\dummyft@=\z@ \fontdimen19 \dummyft@=\z@ \fontdimen20
\dummyft@=\z@ \fontdimen21 \dummyft@=\z@ \fontdimen22
\dummyft@=\z@
\def\fontlist@{\\{\tenrm}\\{\sevenrm}\\{\fiverm}\\{\teni}\\{\seveni}%
 \\{\fivei}\\{\tensy}\\{\sevensy}\\{\fivesy}\\{\tenex}\\{\tenbf}\\{\sevenbf}%
 \\{\fivebf}\\{\tensl}\\{\tenit}\\{\tensmc}}
\def\dodummy@{{\def\\##1{\global\let##1=\dummyft@}\fontlist@}}
\newif\ifsyntax@
\newcount\countxviii@
\def\newtoks@{\alloc@5\toks\toksdef\@cclvi}
\def\nopages@{\output={\setbox\z@=\box\@cclv \deadcycles=\z@}\newtoks@\output}
\def\syntax{\syntax@true\dodummy@\countxviii@=\count18
\loop \ifnum\countxviii@ > \z@ \textfont\countxviii@=\dummyft@
   \scriptfont\countxviii@=\dummyft@ \scriptscriptfont\countxviii@=\dummyft@
     \advance\countxviii@ by-\@ne\repeat
\dummyft@\tracinglostchars=\z@
  \nopages@\frenchspacing\hbadness=\@M}
\def\magstep#1{\ifcase#1 1000\or
 1200\or 1440\or 1728\or 2074\or 2488\or
 \errmessage{\string\magstep\space only works up to 5}\fi\relax}
{\lccode`\2=`\p \lccode`\3=`\t
 \lowercase{\gdef\tru@#123{#1truept}}}

\def\scaletype#1{\mag=#1\relax
 \hsize=\expandafter\tru@\the\hsize
 \vsize=\expandafter\tru@\the\vsize
 \dimen\footins=\expandafter\tru@\the\dimen\footins}

\def\scalefont#1#2\andcallit#3{\edef\font@{\the\font}#1\font#3=
  \fontname\font\space scaled #2\relax\font@}
\def\Mag@#1#2{\ifdim#1<1pt\multiply#1 #2\relax\divide#1 1000 \else
  \ifdim#1<10pt\divide#1 10 \multiply#1 #2\relax\divide#1 100\else
  \divide#1 100 \multiply#1 #2\relax\divide#1 10 \fi\fi}
\def\scalelinespacing#1{\Mag@\baselineskip{#1}\Mag@\lineskip{#1}%
  \Mag@\lineskiplimit{#1}}
\def\wlog#1{\immediate\write-1{#1}}
\catcode`\@=\active

\magnification 1200
\baselineskip 18pt
\def\pbf{\par\bigpagebreak\flushpar}
\def\pmf{\par\medpagebreak\flushpar}
\def\un{\underbar}
\def\suml{\sum\limits}
\TagsOnRight
\def\intl{\int\limits}
\def\qed{\vrule height4pt width3pt depth2pt}
\def\a{\alpha}
\def\b{\beta}
\def\e{\varepsilon}
\def\upn{^{(n)}}


\title
An Unexpected Connection Between Branching Processes and Optimal\\
Stopping
\endtitle

\pbf
\author
David Assaf, Larry Goldstein and Ester Samuel-Cahn\\
Hebrew University, University of Southern California and Hebrew
University
\endauthor

\
\vskip 9cm

\pbf
Abbreviated title:  Branching Processes and Optimal Stopping

\pbf
\un{Key words and phrases}:  Optimal stopping problem, Galton-Watson process, extinction
 probability, varying environment,
inhomogeneous Galton-Watson process, prophet inequalities.

\par\newpage\par\flushpar

\centerline{\smc{Abstract}}

\bigskip

A curious connection exists between the theory of optimal stopping
for independent random variables, and branching processes. In
particular, for the branching process $Z_n$ with offspring
distribution $Y$, there exists a random variable $X$ such that the
probability $P(Z_n=0)$ of extinction of the $n$th generation in
the branching process equals the value obtained by optimally
stopping the sequence $X_1,\ldots,X_n$, where these variables are
i.i.d distributed as $X$. Generalizations to the inhomogeneous and
infinite horizon cases are also considered. This correspondence
furnishes a simple `stopping rule' method for computing various
characteristics of branching processes, including rates of
convergence of the $n^{th}$ generation's extinction probability to
the eventual extinction probability, for the supercritical,
critical and subcritical Galton-Watson process. Examples, bounds,
further generalizations and a connection to classical prophet
inequalities are presented. Throughout, the aim is to show how
this unexpected connection can be used to translate methods from
one area of applied probability to another, rather than to provide
the most general results.

\
\vfill

\vfill
AMS 1991 subject classification: Primary: 60G40, 60J80.

\par\newpage\par\flushpar

\subheading{1. Introduction and Summary}

The purpose of the present note is to highlight what we believe
to be a hitherto unnoticed connection between
two seemingly unrelated
topics in applied probability:
Optimal Stopping Theory for independent random variables, and
Branching Processes and their extinction probabilities.
We show how results in one area can be used to easily establish
results in the other. Our main result is based on a mapping
$Y \rightarrow X$ from integer valued offspring distributions to a
distribution on $[0,1]$ such that the probability of extinction by
generation $n$ of the Galton-Watson branching
process with offspring distribution $Y$ equals the value obtained by optimally stopping a sequence of $n$ independent variables distributed as $X$.
This correspondence is purely analytic, and in particular, we are
not able to present a probabilistic reason, such as a coupling, which
explains it. As the focus is on the `unexplained' connection, in
exploiting the analytic equivalence of the two areas we do not strive
for the most general results, but rather emphasize how one area can
inform another area which is seemingly unrelated.

In Section 2 we outline the basic concepts needed from each of the two
topics. In Section 3 we present our main result, a mapping $Y \rightarrow X$, from
integer valued offspring distributions to distributions on [0,1] such
that the probability
of extinction by generation $n$ of the Galton-Watson branching process
with offspring
distribution $Y$ equals the value obtained by optimally stopping a sequence of $n$
independent variables distributed as $X$.
Examples of this correspondence are given in Section 4.
Section 5 is devoted to proving, by means of ``stopping rule" methods,
various (known) results on rates of convergence of the probabilities of
extinction of the $n^{th}$ generation, denoted
$q_n$, to the eventual probability of extinction, $\pi$, in the subcritical, critical
 and supercritical cases of the Galton-Watson process. In Section 6 we
generalize the results to ``inhomogeneous"
Galton-Watson processes, and provide examples.
In Section 7 we show, in the inhomogeneous case,
 how the use of sub-optimal stopping rules and prophet inequalities
may provide bounds on branching process extinction probabilities, and explore
further connections to the prophet value.
\subheading{2. Basic Concepts}
\pmf
a) \ \un{Optimal Stopping Theory}.  Consider a sequence $X_1, X_2, \dots, X_n$
 of independent random variables with known distributions.
A statistician gets to view the values sequentially, and at each stage must decide
whether to take the present variable or continue. Exactly one variable must be
selected; there is no recall, and hence a variable which has been passed up is no longer available at a later stage to the statistician.
The goal of the statistician is to pick as large a value as possible.
If stopping has not occurred before time $n$ the variable $X_n$ is automatically selected.
The number of variables, $n$, is called the horizon of the problem.
The value to the statistician of using a stopping rule $t$ is
$$
EX_t = E\suml^n_{i=1} X_i I(t=i),
\tag 2.1
$$
where $I$ is the indicator function.
The goal is to maximize the value in (2.1) over all possible stopping rules.

The general theory of optimal stopping is developed in Chow, Robbins and Siegmund
 (1971).
For the finite horizon case an optimal rule always exists and can be obtained
by backward induction.
(See Theorem 3.2, p. 50 of Chow, Robbins and Siegmund (1971)).
In the case of independent random variables the optimal rule has a particularly simple form.
Let $V_i^n$ be the value obtained by optimally stopping the sequence
$X_i,\ldots,X_n$; since stopping must occur at or before time $n$
we set $V^n_{n+1} = - \infty$. If stopping has not occurred by time $i$,
it is optimal to choose $X_i$ only if it is better than or equal to
what is expected in the future.
That is, if $X_i \ge V_{i+1}^n$ the value $X_i$ is selected, and passed up otherwise.
Hence, the value $V_i^n$ is the expectation of the larger of $X_i$ and $V_{i+1}^n$,
that is,
$$
V_i^n= E[X_i \vee V_{i+1}^n].
$$
Alternatively, letting
$$
h_i (a) = E[X_i \vee a]
\tag 2.2
$$
we may write the following recursion for the sequence of values $V_i^n$;
$$
V^n_i = h_i (V^n_{i+1}), \;\;\;\; i = n, n-1, \dots, 1.
\tag 2.3
$$
An optimal stopping rule is
$$
t^*_n = \min \{ i\colon X_i \ge V^n_{i+1}\}.
\tag 2.4
$$
Note that $t^*_n$ will definitely stop by time $n$, if it has not stopped
earlier.
The value of this rule to the statistician is given by $V^n_1$.
In the case where $X_n \ge 0, \; \; V^n_{n+1} = - \infty$ can be replaced by $V^n_{n+1} = 0$.
The case where the $X_i$'s are nonnegative and i.i.d. is of particular interest.
In this case $h_i$ in (2.2) does not depend on $i$, and the index $i$ will
be omitted. Letting
$$
h^{(1)}(a) = h(a) \; \; \hbox{\rm and} \; \; h^{(n+1)}(a) = h (h^{(n)} (a)),
 \;\;\;  n =1, 2, \dots,
\tag 2.5
$$
we have
$$
V_1 ^n = h(V_2^n) = h^{(2)}(V_3^n)=\cdots =h^{(n)} (0).
$$
If we let $V_k$ denote the value for a $k$-horizon problem, then
$V^n_i = V_{n-i+1}, \; \; i = 1, \dots, n$, and
$$
V_k =h^{(k)} (0), \; \quad  k = 1, 2, \dots,.
\tag 2.6
$$
For an infinite horizon problem in this i.i.d. setting, the value
$V_\infty = \lim_{n\to\infty} V_n$ is the supremum over all stopping
rules $t$ with $P(t < \infty)=1$. It equals the rightmost value of the
support of $X$, that is, the essential supermum of $X$.
An optimal rule achieving $V_\infty$ will, however, not exist unless $X$
attains this value with positive probability.

\pbf
(b) \ \un{The Galton-Watson branching process}:

Let ${\cal Y}$ be the set of all nonnegative, nondegenerate integer valued
random variables
excluding the variables for which $P(Y=0)=0$.
For $Y \in {\cal Y}$ let $p_k = P(Y=k), k = 0, 1, \dots$ and
$$
g(s) = \suml^\infty_{k=0} p_k s^k
\tag 2.7
$$
be the generating function of $Y$, which is well defined for $0 \le s \le 1 $,
with $g(0) = p_0$ and $g(1) = 1$.
Note that if $EY<\infty$ then $g'(1) = EY$, and if $EY^2<\infty$ then
$g''(1) = EY^2 -EY$.
All derivatives of $g(s)$ for $s \in [0, 1)$ exist and are nonnegative,
thus in particular $g(s)$ is increasing and convex; the function $g$ will be
strictly convex unless it is linear, that is, unless $p_0+p_1=1$.

For given $Y_n \in {\cal Y}$, $n=1,2,\ldots$, define the
(inhomogeneous, or varying environments)
Galton-Watson branching process, with offspring distribution $Y_n$ at generation $n$,
as the discrete time stochastic process $\{ Z_n\}^\infty_{n=0}$
with $Z_0 = 1$ and
$$
Z_{n+1} = \suml^{Z_n}_{i=1} W_{ni},
\tag 2.8
$$
where $W_{ni}$ are i.i.d. distributed like $Y_n$. The value $Z_n$ is the
size of the $n^{th}$ generation of a population which
begins with a single individual at time 0, where each member
of generation $n$ gives rise to offspring for the next generation with
distribution $Y_n$, independently of all the other members.
Letting $g^{(n)}$ be the generating function of $Z_n$, and $g_n$ be the
generating function of $Y_n$, we have the well known relation
$$
g^{(n)} (s) = g_1(g_2(\cdots g_n(s))),
\tag  2.9
$$
which can be verified by induction.
A quantity of major interest is the probability that the $n^{th}$ generation is extinct
$$
P(Z_n=0)=g^{(n)}(0) = q_n.
\tag 2.10
$$
Since $Z_n = 0$ implies $Z_{n+1} = 0$,
we have $0 \le \tilde q_1 \le \tilde q_2 \le \dots $, and thus
$\lim_{n \rightarrow \infty} \tilde q_n = \tilde \pi \le 1$
exists. The limit $\tilde \pi$ is the probability of eventual extinction.
Furthermore, it is easily seen that
$EZ_n$, the expected size of generation $n$,
equals $\prod_{j=1}^n EY_j$. This follows
by computing $[g^{(n)}]'(1)$ in (2.9) and using $g_j(1)=1$.

When all the $Y_n$ have identical distributions with generating function $g$,
$$
g^{(1)} (s) = g(s), \;\;  g^{(n+1)}(s) = g(g^{(n)}(s)) \;\;\;
n=1,2,\ldots,
\tag 2.11
$$
and we denote $q_n=P(Z_n=0)=g^{(n)}(0)$, and $\lim_{n \rightarrow \infty}q_n =\pi$.
As is well known (see e.g. Karlin and Taylor (1975), Chapter 8) in this instance
$\pi$ is the smallest root of the equation
$$
g(s) = s.
\tag 2.12
$$
The value $s=1$ is always a root of (2.12), and
it is the smallest root if and only if $EY<1$ (the subcritical  case),
 or $EY = 1$ (the critical case).
There is positive probability of never becoming extinct,
that is, of having $\pi < 1 $, iff $EY>1$ (the supercritical case).

\subheading{3. Connection}
Our main result in the present section is to exhibit the connection between Optimal Stopping
and homogeneous Galton-Watson Processes. In particular we link the optimal stopping value
$V_n$ to the extinction probability $q_n$ using the following Theorem.
\proclaim{Theorem 3.1}
Let $Y \in {\cal Y}$ have generating function $g$, and let $\pi$ be the smallest root of
the equation $g(s)=s$. Then the function $F(x)$ given by
$$
F(x)=
\cases
0 &x<0\\
g'(x) &0\le x < \pi\\
1 &\pi \le x,
\endcases
\tag 3.1
$$
is a distribution function. Let $X$ have distribution (3.1), and
$h(a) = E[X \vee a]$. Then
$$
h(a) = g(a) \; \; \hbox{\rm for} \; \; 0 \le a \le \pi.
\tag 3.2
$$
Also
$$
EX=P(Y=0), \quad \;\; P(X=0)=P(Y=1) \quad \;\;
\hbox{and when} \; \;\pi=1,\quad P(X < 1)=EY.
\tag 3.3
$$
The variable $X$ has an atom of size $P(Y=1)$
at 0, an atom of size $1-g'(\pi)$ at $\pi$, and density $g''(x)$ on  $(0, \pi)$.
Exactly one distribution satisfies (3.2).
\endproclaim
\demo{Proof}
The function $g'(x)$ is non-negative and nondecreasing
for $0 \le x \le \pi$. Note that $g' (\pi) \le 1$, since $g$ is convex and $s < g(s)$ for
all $0 \le s < \pi$. Further, by definition, for $0 \le a \le \pi$,
$$
h(a) = E[X\vee a]=\intl^\infty_0 P([X \vee a]>x) dx = \pi - \intl^\pi_a g'(x) dx =
\pi - g(\pi) + g(a) = g(a),$$
which is (3.2).
For $a=0$, (3.2) yields $h(0)=g(0)$, or $EX=P(Y=0)$, and
$P(X=0)=g'(0)=P(Y=1)$. When $\pi=1$, $EY \le 1$
and (3.1) yields $EY = g'(1) = F(1^-)=P(X <1)$.

To show uniqueness, suppose (3.2) holds for some $X^*$ with distribution function $F^*$.
Since $\pi = g(\pi) = h(\pi) = E[X^* \vee \pi]$, it follows that $P(X^* > \pi) = 0$,
 i.e. $F^* (x) = 1$ for all $x \ge \pi$.
Also, since $g$ is differentiable in $0< s < \pi$, so is $h$.
But for  $0 \le s \le \pi,  h(s) = E[X^*\vee s] = 1-\intl^1_sF^*(x) dx = g(s)$,
 thus $g'(s) = F^* (s)$ for $0<s<\pi$, and thus, by right continuity, $F^*(x) = F(x)$
for all $x$.
\hfill\qed

\proclaim{Theorem 3.2}
Let $Z_n$ be a Galton-Watson process with offspring distribution $Y \in {\cal Y}$
and extinction probability $q_n = P(Z_n=0)$.  Let
$X_1,\ldots, X_n$ be i.i.d. with distribution function (3.1), and let $V_n$ be its
optimal stopping value.
Then,
$$
V_n = q_n, \; \; n =1, 2, \dots
\tag 3.4
$$
\endproclaim
\demo{Proof}
For $0 \le a \le \pi$ we have $0 \le g (a) \le \pi$. By (3.2) and
induction,
$$
h^{(n)} (a) = g^{(n)} (a) \; \; \hbox{\rm for } \; \; 0 \le a \le \pi.
\tag 3.5
$$
Using (2.6) and (2.10) and setting $ a = 0 $ in (3.5) yields (3.4).
\hfill\qed
\demo{Remarks}

{\bf 3.1} Equality (3.2) cannot hold for $\pi<a<1$ since in
this interval $g(a)<a$, while $h(a)=E[X \vee a] \ge a$.

{\bf 3.2} The distribution of $Y$ is uniquely determined by the sequence $\{ q_n\}^\infty_1$,
since an analytic function $g$ is uniquely determined by its values on an infinite
sequence of values having a limit point.
Thus there are no two different $Y$'s with the same $q_n$-sequence.

{\bf 3.3} In contrast to Remark 3.2, there are many different i.i.d. sequences
 of $X$'s with  values $\{ V_n\}^\infty_1$.
For a construction, see Hill and Kertz (1982).

{\bf 3.4} We excluded from ${\cal Y}$ the variables for which $P(Y=0) = 0$.
For such variables $\pi = 0$ is the smallest root of (2.12).
Note that for this case $F$ of (3.1) gives unit mass to $0$, thus (3.2) and (3.4)
are formally true also for this case.

{\bf 3.5}  Theorem 3.1 shows that for each $Y \in {\cal Y}$
there exists an $X$ taking values in [0,1] such that (3.2) holds. However, it is
not true that for each $X$ taking values in [0,1] there exists a corresponding
$Y \in {\cal Y}$. Necessary and sufficient conditions for $X$ to
correspond to a $Y \in {\cal Y}$ is that $X$ has a distribution function
$F$ of the form
$$
F(x) = \cases 0 &x<0\\
k(x) &0\le x < \pi\\
1 &\pi\le x, \endcases \tag 3.6
$$
for some $0<\pi \le 1$,
and that
\pmf
(i) \ $k( \; )$ has a power series expansion with all coefficients nonnegative, and
\ (ii) \ There exists a constant $c > 0$ such that $g(s) = \intl^s_0 k(x)
dx + c$
satisfies (a) $g(1) = 1$, \ (b) \ $g(\pi) = \pi$. This
fact suggests that it will be easier to use the correspondence to translate
properties of optimal stopping into properties about Galton-Watson
processes,
than vice versa.

\subheading{4. Examples}

The correspondence between $Y$ and $X$ of (3.1), yields some interesting
relationships.

\demo{Example 4.1}$Y \sim {\cal B}(p)$ Bernoulli. In this case
$P(Y=1) = p = 1 - P(Y=0)$ and
clearly $\pi=1$. As $g(s) = (1-p) + ps$, $F(s) = p$ for $0 \le s <1$,
and $F(1)=1$. Hence, $X \sim {\cal B}(1-p)$.

\demo{Example 4.2}  $Y \sim m{\cal B}(p)$, $m \ge 2$, that is,
$P(Y=m) = p = 1-P(Y=0)$, $g(s) = (1-p) + ps^m$, and $EY=mp$.
Using (3.3), since $P(Y=1)=0$, $X$ has no mass at zero,
but has mass $1-g'(\pi)=1-mp\pi^{m-1}$ at $\pi$. Therefore, for $0 \le s \le  \pi$,
$$
F(s)=mp\pi^{m-1} \left( \frac{s}{\pi}\right)^{m-1} +
I(s=\pi) (1- mp \pi^{m-1}),
$$
that is, $X$ is a mixture of $\mathop{\max} \limits_{i=1, \cdots, m-1} U_i$, where
$U_i$ are i.i.d. $U(0,\pi)$, with probability $mp \pi^{m-1}$, and a point
mass at $\pi$ with probability $1-mp \pi^{m-1}$.
In particular, for $m = 2$ (corresponding to a splitting of a cell), in
the critical case $p=1/2$, $X \sim U (0,1)$.
For $m = 2$ and the supercritical case $ p > 1/2$, the eventual extinction
probability is the smallest solution to $1-p + ps^2 -s=0$, which is $\pi = (1-p)/p$.
Therefore, $X$ is a mixture of  $U(0,\pi)$ variable with probability $2(1-p)$
and a point mass at $\pi=(1-p)/p$ with probability $2p-1$. In the subcritical case $p < 1/2$,
$X$ is a mixture of a uniform $U(0,1)$
 variable with probability $2p$,
and point mass at 1 with probability $1-2p$.
\demo{Example 4.3}  $Y \sim {\cal P}(\lambda)$, $Y$ is Poisson with parameter $\lambda$,
and $g(s) = e^{\lambda(s-1)}$.
For $\lambda > 1$, $\pi<1$ is the smallest root of $e^{\lambda (s-1)}=s$;
for $\lambda \le 1 $, $\pi=1$. The distribution function of $X$ is
$$
F(x) = \cases 0 &x<0\\
\lambda e^{\lambda (x-1)} &0\le x < \pi\\
1 &x \ge \pi. \endcases
$$
\demo{Example 4.4}  $Y \sim {\cal GG}(b,c)$, Generalized Geometric distribution:
$P(Y=k) = bc^{k-1}, k = 1, 2, \dots $ and $P(Y=0) = 1-\suml^\infty_{k=1} bc^{k-1}
 = (1-b-c)/(1-c)$, for any $b,c>0$ such that $b+c<1$.
The standard geometric
distribution ${\cal G}(p)$ with success probability $p \in (0,1)$, $p+q=1$, is the special
case ${\cal GG}(pq,q)$.
Here $g(s) = P(Y=0) + bs/(1-cs)$ and can be written as
$$
g(s) = (\a + \b s)/(\gamma + \delta s)
\tag 4.1
$$
with
$$
\a = 1-(b+c),\; \;  \b = b-c(1-c), \; \; \gamma = 1-c, \; \; \delta = - c(1-c).
\tag 4.2
$$
This is (according to Athreya and Ney (1972, p. 6))
essentially the only nontrivial example where $g^{(n)}(s)$, and hence
$g^{(n)} (0) = q_n$, can be computed explicitly.
This example is also discussed  in most other texts on branching processes, see e.g.
 Harris (1963, p. 9), and Karlin and Taylor (1975, p. 402).
See also the continuation of this example in Example 6.2, below.
Since $EY=b/(1-c)^2$ it follows  easily that for $b>(1-c)^2$ the eventual extinction
 probability is $\pi = [1-(b+c)]/c(1-c) = -\a/\delta$.
In all other cases  $\pi = 1$.
Here $X$ has c.d.f.
$$
F(x) = \cases 0 &x<0\\
b/(1-cx)^2 &0 \le x < \pi\\
1 &x  \ge \pi.\endcases
\tag 4.3
$$
\subheading{5. Convergence rates of the extinction probabilities for the Galton-Watson
process}

The purpose of the present section is not to derive new results,
 but rather to show how well-known results in branching theory have simple proofs by
means of stopping rules.
We do not strive for the most far-reaching results, and are content with obtaining
 rates for which $q_n \to \pi$.
\proclaim{Theorem 5.1}
\pmf
(a) Supercritical case: If $EY>1$ then  $\pi<1$ and
$$
0<\pi-q_n < \pi [g'(\pi)]^n.
\tag 5.1
$$
\pmf
(b) Subcritical case: If $EY<1$ (and $P(Y=0)<1)$,
then $\pi=1$ and
$$
0<1-q_n \le [EY]^n,
\tag 5.2
$$
and the inequality on the right in (5.2) is strict if and only if $P(Y\le 1) < 1$.
\pmf
(c) Critical case: If $EY=1$, Var($Y) = \sigma^2 <\infty$ then
$$
\lim\limits_{n\to\infty} n[1-g'(q_n)]=2,
\tag 5.3
$$
or equivalently,
$$
\lim\limits_{n\to\infty} n(1-q_n) = 2/\sigma^2.
\tag 5.4
$$
More generally, if $EY=1$ and
$$
\lim\limits_{s\to 1^{-}} (1-s)g''(s)/[1-g'(s)]=\a
\tag 5.5
$$
for some $0<\a\le 1$, then
$$
\lim\limits_{n\to\infty} n[1-g'(q_n)] = 1+\a^{-1}.
\tag 5.6
$$
\pmf
\endproclaim
\demo{Proof}  (a)  According to Theorem 3.1, for $X$ corresponding to $Y$,
$P(X=\pi) = 1-g'(\pi)$, which is positive.
Now consider the suboptimal stopping rule $t$ which stops at the smallest $i$
 for which $X_i = \pi$, and if no such $i$ exists, stops at time $n$
 anyway. Since this rule is suboptimal, $EX_t$, the expected value to the statistician
using rule $t$, is at most $V_n$,
 but is greater than $\pi$ times the probability that the value $\pi$ will be observed,
 since stopping at $t = n$ with some value smaller than $\pi$ will still yield
 a positive expected return.
The probability of never observing a value $\pi$ is $[g'(\pi)]^n$.
Thus $\pi (1-[g'(\pi)]^n)< EX_t \le V_n=q_n$,
 from which (5.1) follows.
\pmf
(b) The proof of (b) is essentially the same as (a), using
$\pi = 1$, and $P(X=1) = 1-g'(1) = 1 - EY$.
Equality in (5.2) holds if and only if the ``suboptimal" rule $t$
is actually optimal.
This happens if and only if $X$ is Bernoulli.
This case is described in Example 4.1, where $P(Y \le 1) = 1$, and by the
uniqueness of $X$, as stated in Theorem 3.1, this is the
only case.
\pmf
(c) We shall draw on the results of Kennedy and Kertz (1991), who show that
 the asymptotic behavior of the value sequence $V_n$ for optimal stopping of i.i.d. random variables depends on to which extremal distribution domain
$X$ belongs. In the present case, $X$ has no mass at 1, is bounded above
by 1, has distribution function $g'(x)$, and the non-zero
density $g''(x)$ for $0 < x < 1$. In terms of the given c.d.f. and
density, condition (5.5) is equivalent to
the condition for a Type III extreme value distribution given in Theorem 1.6.1. of
Leadbetter, Lindgren and Rootz\'en (1983)  (See also e.g. de Haan (1976),
Theorem 4 and the remark which follows).
Theorem 1.1 of Kennedy  and Kertz (1991) now yields (5.6).
Note that when Var($Y) = \sigma^2 < \infty$ then $g''(1) = \sigma^2$,
 (since $EY =1)$, and the value of the limit in (5.5) is necessarily 1.
Thus (5.3) is the particular case of (5.6) with $\a = 1$.
Note that by convexity the value in the left hand side of (5.5)
 for every fixed $s$ is necessarily less than 1, and hence only $\a$-values less than
or equal to one can be obtained as limits in (5.5).

To see that (5.3) is equivalent to (5.4), note that since $q_n \to 1$ and
$\lim_{n \rightarrow \infty} (1-g'(q_n))/(1-q_n) =
\sigma^2$ we obtain $\lim_{n \to \infty}
n(1-g'(q_n)) = \lim_{n \to \infty}  n (1-q_n)\sigma^2$.
\hfill \qed
\demo{Remarks}

{\bf 5.1} Standard proofs of various parts of Theorem 5.1 can be found in most standard
 texts in Branching processes.

{\bf 5.2} We see that the convergence of $q_n$ to $\pi$ is at a geometric rate in both
 the supercritical and subcritical cases.
It is at the order of $0(1/n)$ in the critical case
when Var($Y) < \infty$, but $1-q_n$ converges to zero faster when Var($Y) = \infty$.

{\bf 5.3}  The branching process with $EY = 1$ and Var($Y)=\infty$ is studied in Slack, (1968).
Note that all values of $\a, 0 < \a \le 1 $ can be attained as the
limit
 in (5.5), as seen from the following
\demo{Example 5.1}
For $0 < \a \le 1$, and $0<c \le 1/(1+\a)$, let $Y$ have generating function
$$
g(s) = s+(1-s)^{1+\a} c.
\tag 5.7
$$
It is easily seen that this corresponds to the distribution
$$
\eqalign{
&P(Y=0) = c, \; P(Y=1) = 1- (1+\a)c\cr
&P(Y=k) = (-1)^k c \prod_{j=-1}^{k-2} (\a-j)/k!, \; \; k = 2, 3, \dots\cr}
\tag 5.8
$$
(For $\a=1$ it follows that $P(Y=k) = 0 $ for $k>2$).
Since $g'(1)=1$ it follows that $EY=1$ and easy arithmetic yields (5.5).
Here (5.6) can be stated as
$$
\lim\limits_{n\to\infty} n(1-q_n)^\a = (c\a)^{-1}
\tag 5.9
$$
and shows that $q_n$ tends to 1 faster, the smaller $\a$.
(Note that for $\a=1$ one has Var($Y) = 2c$ and (5.9)
agrees with (5.4) in this case).

{\bf 5.4}  Though in most natural situations the limit in (5.5) does exist,
one can exhibit generating functions for
which the limit in (5.5) fails to exist. One such construction is a
function having a
coefficient sequence which essentially alternates between the
coefficient sequences of generating
functions of the form (5.7) for two different values of $\a$.

\subheading{6. Inhomogeneous branching processes}

In this section we consider the inhomogeneous
branching process, as presented in Section 2(b). Here the
offspring distribution in generation $i$ is $Y_i$, where the $Y_i$ need not
have identical distributions.
To each $Y_i$ there is a corresponding $X_i$
defined through (3.1), where $g$ there is replaced by $g_i$, and $\pi$ by $\pi_i$
(where $\pi_i$ is the eventual extinction probability of an ordinary
Galton-Watson
process with fixed offspring distribution $Y_i$.)
Now consider an optimal stopping problem where $X_1, \dots, X_n$ are observed
sequentially.
From (2.2) and (2.3) it follows that the value $V^n_1$ to the statistician,
 of this sequence is
$$
V^n_1 = h_1 (h_2(\cdots h_n(0))).
\tag 6.1
$$
If we denote more generally
$$
h^{(n)}(a) = h_1 (h_2(\cdots h_n(a)))
\tag 6.2
$$
then, using (2.9), we can generalize Theorem 3.1 and (3.5) as follows.
\proclaim{Theorem 6.1} Suppose $\pi_1 \ge \pi_2 \ge \cdots \ge \pi_n$.
Then
$$
h^{(n)} (a) = g^{(n)} (a) \; \; \hbox{\rm for} \; \; 0 \le a \le \pi_n,
$$
and thus also $V^n_1 = \tilde q_n$.
\endproclaim
The proof is straightforward and hence omitted.

Inhomogeneous Galton-Watson processes have been studied quite extensively in the
literature.
The earlier references are Jagers (1974) and Jirina (1976).
See also  Section 3.5 in Jagers (1975).
One of the latest references we have come across  is D'Souza (1995).
See also all related references mentioned there.
All papers deal with various aspects of the limiting value of $Z_n$ under different
 assumptions on the $Y_i$'s.
The following theorem has a very simple ``stopping rule" proof.
\proclaim{Theorem 6.2}  Suppose $\pi_i =\pi_0$ for all
$i = 1, 2, \cdots$, and denote $r_i = P(Y_i \not =1)$.
Then $\tilde \pi \le \pi_0$ with
$$
\tilde \pi = \pi_0\; \;  \hbox{\rm if and only if}
\; \; \suml^\infty_{ i=1} r_i = \infty.
\tag 6.3
$$
\endproclaim
\demo{Proof}
Let $X_i$ correspond to $Y_i$ through the relation (3.1). By
Theorem 6.1, $\tilde q_n=V_1^n$, and hence $V_\infty =
\mathop{\lim}\limits_{n\to\infty} V^n_1 =
\mathop{\lim}\limits_{n\to\infty} \tilde q_n = \tilde \pi$.
Since by (3.1) $X_i \le \pi_0$, $i=1,2,\ldots,$, $\tilde \pi=V_\infty \le \pi_0$.

We will show that $\tilde \pi = \pi_0$ if and only
if for all $0 < \e < \pi_0$,
$$
\suml^\infty_{i=1} [1-g'_i (\pi_0-\e)]=\infty,
\tag 6.4
$$
and then prove that (6.4) is equivalent to the condition
$\sum r_i = \infty$. Note that
$P(X_i \ge \pi_0 - \e) = 1-g'_i(\pi_0-\e)$.
Thus if (6.4) holds then
$P(X_i \ge \pi_0 - \e \; \hbox{infinitely often}) = 1$.
Hence, for the rule $t=\inf\{i\colon X_i\ge \pi_0-\e\}$, we have
$P(t < \infty)=1$, and
the value for this rule, $EX_t$, is at least $\pi_0-\e$, and hence
$\pi_0 -\e \le V_\infty \le \pi_0$.
Since this is true for every $\e>0$ it follows that
$\tilde \pi=V_\infty =\pi_0$.
Conversely, if (6.4) fails for some $0 < \e_0 < \pi_0$ then by the Borel-Cantelli lemma,
$P(X_i \ge \pi_0-\e_0  \; \hbox{infinitely often})<1$,
and hence there is positive probability of
never seeing a value greater than $\pi_0 - \e$, thus
the supremum of the expected return over
all stopping rules is less than $\pi_0$, that is,
$\tilde \pi < \pi_0$.

It remains to verify that (6.4) is equivalent to the condition
$\sum r_i = \infty$.
Let $P(Y_i = k) = p_{ik}, \; \; k = 2, 3, \ldots $.
Then
$$
g'_i (s) = 1-r_i + \suml_{k\ge 2} kp_{ik} s^{k-1}.
$$
Note that for any $k\ge 1$, \  $ (\pi_0 - \e)^k \le (\pi_0 -\e)\pi_0^{k-1} = (1-\e/\pi_0)\pi_0^k$,
and also that $g_i'(\pi_0) \le 1$; thus
$$
\eqalign{
&1-g'_i(\pi_0-\e)=r_i-\suml_{k\ge 2} kp_{ik} (\pi_0 - \e)^{k-1} \ge r_i - (1-\e/\pi_0)
\suml_{k\ge 2} kp_{ik} \pi^{k-1}_0\cr
&=r_i-(1-\e/\pi_0) [g'_i(\pi_0) -1+r_i]
\ge r_i-(1-\e/\pi_0)r_i=(\e/\pi_0)r_i.\cr}
\tag 6.5
$$
Thus, if $\sum r_i = \infty $, (6.4) holds for every $0< \e < \pi_0 $
 (and $\tilde \pi = \pi_0$).
On the other hand, from the first equality in (6.5) it follows that
$1-g'_i (\pi_0 - \e) \le r_i$.
Thus $\sum r_i < \infty$ implies that the sum in (6.4) converges.
\hfill \qed

\demo{Remark 6.1} Suppose $EY_i \le 1 $ for all $i$, and $\suml^\infty_{i=1} [1-EY_i]
=\infty$.
Since $1-EY_i \le r_i$,  $\suml^\infty_{i=1} r_i = \infty$ and $\tilde \pi=1$ follows.

Since $EY_i \le 1$, one has
$P(X_i=1)= 1- g'(1) = 1 - EY_i$.  Thus
the probability that $X_i = 1$ infinitely often, equals one,
so the rule which stops for the smallest
$i$ for which $X_i = 1$, stops with probability 1. Thus the value $V_\infty = 1$ is,
 in this case, attainable by a stopping rule $t$ with $P(t < \infty) =1$.
In all other situations where $\sum r_i = \infty$, the value $V_\infty=1$ is not
attainable, and only $\e$-optimal stopping rules exist.

\demo{Remark 6.2}  Note that $\pi_i = \pi_0$ for all $i$
implies by (2.9) that $g^{(n)} (\pi_0) =\pi_0$ also.
Thus unlike the situation in the homogeneous Galton-Watson process
where $\lim\limits_{n\to\infty}
 g^{(n)} (s) = \pi $ for $0 \le s \le \pi$, in the inhomogeneous case,
it may happen that even though $\lim\limits_{n\to\infty} g^{(n)} (s)
$ exists for all $0 \le s \le 1$, this limit need not equal $\pi_0$ for
$0 \le s < \pi_0$, unless $\sum r_i = \infty$.
A similar remark is true also for the case $\pi_0 = 1$.
\demo{Example 6.1} Let $Y_i$ take the values $0, 1 $ and $2$ only, with probabilities
 $P(Y_i=0) = r_i/3$, $P(Y_i=1) = 1-r_i$ and $P(Y_i = 2) = 2r_i/3$, where $0 < r_i \le 1$.
Here $g_i (s) = r_i/3+(1-r_i)s+2r_is^2/3$, and it is easily checked
 that $\pi_i = 1/2$ for all $i$.
Note that here $EY_i = 1+r_i/3$.
Since $EZ_n = \prod\limits^n_{i=1} EY_i$, the condition $\sum r_i = \infty$ is equivalent
 to $\lim\limits_{n\to\infty} EZ_n = \infty$.
\demo{Example 6.2} As in Keiding and Nielsen (1975), let $Y_i$
 have the Generalized Geometric distribution, ${\cal GG}(b_i,c_i)$,
 as described in Example 4.4; hence, $Y_i$ has generating function as in (4.1)
$$
g_i(s) = (\a_i +\b_i s) /(\gamma_i + \delta_is)
\tag 6.6
$$
where the constants are defined as in (4.2).
Then it can be verified by induction that
$$
g^{(n)} (s) = (\a^{(n)} +\b^{(n)}s)/(\gamma^{(n)}+\delta^{(n)}s),
\tag 6.7
$$
and the values of $\a^{(n)}, \b^{(n)}, \gamma^{(n)}$ and $\delta^{(n)}$
can be obtained explicitly.
We shall consider in detail the case where all $Y_i$ are ``critical",
i.e. $b_i = (1-c_i)^2$.
For this case  let
$$
S^{(n)}_i = \sum c_{k_1} \cdots c_{k_i}
\tag 6.8
$$
where the summation is over all $1 \le k_1 < \cdots < k_i \le n$.
Set $S^{(n)}_0 = 1$.
Then one can verify that
$$
\eqalign{
&\a^{(n)} = \suml^n_{j=1} (-1)^{j-1}jS^{(n)}_j, \qquad  \b^{(n)}=\suml^n_{j=0}
 (-1)^j (j+1)S^{(n)}_j\cr
&\gamma^{(n)}=1+\suml^n_{j=2} (-1)^{j-1} (j-1)S^{(n)}_j, \qquad
\delta^{(n)}=-\a^{(n)}.\cr}
\tag 6.9
$$
Clearly here $\tilde q_n = \a\upn/\gamma\upn$.
It follows from (6.8) and (6.9) that for any $n$, all $n!$ permutations of the order
 of the $Y_i$s
 yield the same distribution for the $n^{th}$ generation, $Z_n$.
Note that here $P(Y_i=1) = 1-r_i = (1-c_i)^2$ which implies that $r_i = c_i
 (2-c_i)$.
Thus, by Theorem 6.2, $\tilde \pi=1$ if and only if $\sum c_i = \infty$.

It is of interest to note that the permutation invariance mentioned above can
 be generalized.
Let $Y_1 $ and $Y_2$ have generating function of the form (4.1).
Then  the  generating function of $Z_2$ is (see (2.9))
$$
g_1 (g_2(s)) \, = \, {(\a_1\gamma_2+\b_1\a_2)+(\a_1\delta_2+\b_1\b_2)s\over
(\gamma_1\gamma_2+\delta_1\a_2)+(\gamma_1\delta_2+\delta_1\b_2)s},
$$
and it can then be verified that $g_1(g_2(s)) = g_2(g_1(s))$ if and only if $\a_1/\delta_1
=\a_2/\delta_2$.
But for $\pi_i<1$ one has
$-\a_i/\delta_i=\pi_i$, thus the order does not matter if and only if
$\pi_1=\pi_2$.
This generalizes immediately for composing $n$ such generating functions,
 and shows that the order
 of the $Y_i$s does not matter if and only if all $\pi_i=\pi_0 < 1$ in this case.
We do not know if this property has been observed earlier.
Translating to optimal stopping, we have obtained a sequence of non identically distributed
variables for which the optimal stopping value is the same, no matter
in which order the variables appear.

It is easy to show, by working out the distribution of $Z_2$ in Example 6.1,
that even though $\pi_1=\pi_2=1/2$ there, the $Y_i$s there do not have
 the permutation invariance property.

\subheading{7. Connections to Prophet Values and Prophet Inequalities}

When not all $\pi_i$ are equal,
or when the necessary condition of Theorem 6.2 fails, one may still obtain
meaningful, though sometimes crude, lower and upper bounds on $\tilde \pi$
through the use of suboptimal stopping rules, the `prophet' value
and the `prophet inequality.' If $EX_t$ is the value of
any (optimal or suboptimal) stopping rule $t$,
for the $n$-horizon case, then
$EX_t \le V_n=\tilde q_n \le V_\infty = \tilde \pi$, and if $EX_t$ is
the value of a suboptimal rule for
the infinite horizon case, $EX_t \le V_\infty = \tilde \pi$, yielding lower bounds
on $\tilde q_n$ and $\tilde \pi$.
Let $V_p^n = E(\max
 (X_1, \dots, X_n))$  and $V_p^\infty = \lim_{n \rightarrow \infty} V_p^n$.
$V^n_p$ and $V^\infty_p$ are called ``prophet
values". The term ``prophet value" stems from the fact that an individual with
complete foresight of the future would simply select the largest $X_i$ value
in the sequence, and obtain the expected return $V_p$, the ``prophet
value". The prophet values $V_p^n$ and $V_p^\infty$
are usually much easier to compute than the optimal stopping value.
Since the value of any stopping rule is necessarily less than or equal
to that of the prophet, we have the upper bound
$\tilde q_n = V_n \le V_p^n \le V_p^\infty$ and hence ${\tilde \pi} \le
V_p^\infty$.
In addition, the prophet value can also be used to obtain a lower bound on $\tilde \pi$.
It is well-known, (see e.g. Hill and Kertz (1981)) that for a sequence of
nonnegative independent random variables $V_p^n < 2 V_n$, and thus
$V_p^n/2 < V_n = \tilde q_n \le \tilde \pi$
serves as a lower bound on $\tilde q_n$ and $\tilde \pi$. Letting $n
\rightarrow \infty$, we see that
$V_p^\infty /2$ is also a lower bound on $\tilde \pi$.

\demo{Example 7.1} Consider Example 6.1 with $\sum r_i <\infty$.
Since
$g_i'(0)=1-r_i$, and $\pi_i=1/2$, the $X_i$ corresponding to $Y_i$ has mass
$1-r_i$ at zero and is bounded above by 1/2. Hence, the variable $X_i^*$ where
$P(X^*_i = 0)=1-r_i$ and $P(X^*_i=1/2)
= r_i$ is stochastically larger than $X_i$, and therefore
the prophet value for the  $X_i^*$ sequence is an upper bound
on the prophet value for the $X_i$ sequence.
The prophet value for the $X^*_i$-sequence is 1/2 the probability
that any of the $X_i^*$ variables equals 1/2, i.e.,
$({1\over 2})[1-\prod^\infty_{i=1}
 (1-r_i)]$.
To obtain a lower bound on $\tilde \pi$, consider the suboptimal rule which
stops for the
smallest $i$ such that $X_i>0$.
It should be noted that since $\sum r_i < \infty$, this rule does not stop with
probability one unless $r_i=1$
 for some $i$.
Even if $r_i < 1 $ for all $i$
the value of this ``rule" equals the limit of the value of the rule $t_n$ which
 stops for the smallest $i$ such that $X_i >0$, and stops at time $n$
 if no positive $X_i$ is observed up to and including time $n$.
The conditional expected return for stopping at $X_i$, given $X_i>0$, is $1/3$.
Thus the value of this rule is $({1\over 3})[1-\prod^\infty_{i=1} (1-r_i)]$.
A different lower bound can be obtained through the rule which stops for the smallest
$i
$ such that $X_i = 1/2$, if such an $i$ exists.
Its expected return is $({1\over 2})[1-\prod^\infty_{i=1} (1-r_i/3)]$.
Thus
$$
\max\{ \frac{1}{2} [1-\mathop{\Pi}^\infty_{i=1} (1-r_i/3)], \frac{1}{3} [1-\mathop{\Pi}^\infty_{i=1}
(1-r_i)]\} \le
\tilde \pi < ({1/2})[1-\mathop{\Pi}^\infty_{i=1} (1-r_i)].
$$
For example, if $r_i = 1/(i+1)^2$ we have
$$\prod^\infty_{i=1} (1-r_i) = \lim\limits_{n\to\infty}
 \prod^n_{i=1} {i(i+2)\over (i+1)^2} =
\lim\limits_{n\to\infty} \frac{(n+2)}{2(n+1)} = 1/2,$$
so that $1/6\le \tilde \pi\le 1/4$. (Recall that $\pi_0 = 1/2$).

We have shown how the correspondence between $Y$ and $X$ can be used to
obtain information about branching processes from computations involving an
optimal stopping problem. The following theorem shows how the correspondence
can be applied in the other direction.

\proclaim{Theorem 7.1}  Let $Y\in {\cal Y}$ with $EY \le 1$,
and let $X$ be the corresponding random variable with distribution given in (3.1).
With $X_1, \dots, X_n$ i.i.d. random variables distributed like $X$, let $X^*_n
 = \max (X_1, \dots, X_n)$. Then $X^*_n$ corresponds to a $Y^*_n \in {\cal Y}$,
and the prophet value $E X^*_n$ can be computed using
$$
EX^*_n = P(Y^*_n = 0).
\tag{7.1}
$$
\endproclaim
\demo{Proof}  The distribution function of $X^*_n=\max (X_1, \dots, X_n)$ is
$$
F^*_n(x) = \cases
0&x<0\\
[g'(x)]^n &0\le x < 1\\
1&1\le x. \endcases
\tag 7.2
$$
Clearly $k(x) = [g'(x) ]^n$ satisfies condition (i) of Remark 3.5.
Now since $g'(x) \le g'(1)=EY \le 1$, $[g'(x)]^n \le g'(x)$ for $0<x<1$ and
$\intl^1_0 [g'(x)]^n dx\le \intl^1_0 g'(x) dx < 1$. Hence
$$
g^*(s)= \int_0^s k(x)dx +(1-\int_0^1 k(x)dx)
\tag 7.3
$$ further satisfies $g^*(0) > 0$ and
$g^*(1)=1$. Since here $\pi = 1$, condition (ii) (b) of Remark 3.5, $g(\pi)=\pi$, is
equivalent to (ii) (a). \hfill \qed

\demo{Remark 7.1} If $EY>1$ i.e. $\pi<1$, then $\max (X_1, \dots, X_n)$ does not correspond
 to any $Y^*\in {\cal Y}$ since the distribution corresponding to (7.2) for this case
 cannot satisfy (ii) (a) and (b) of Remark 3.5 simultaneously.
\demo{Remark 7.2} When $EY\le 1$, then $EY^*_n =dg^*(x)/dx|_1=k(1)=[g'(1)]^n= [EY]^n$.

For the cases below which illustrate Theorem 7.1, the given $g(x)$ is
sufficiently unaltered upon differentiation, taking powers, and integration
that $g^*(x)$
of (7.3) correspondonds to a variable $Y_n^*$ of the same `type' as
the original $Y$, with a mass at zero according to the constant term in
(7.3).

\demo{Example 7.2} Let $X$ be the variable corresponding to the $Y$
of Example 4.2 with $p\le 1/m$, where $g'(x) =mpx^{m-1}$. Hence
$k(x)=(mpx^{m-1})^n$ for $0 \le x < 1$ and hence
$$
g^*(x) = (mp)^n x^{n(m-1)+1}/[n(m-1)+1]+(1-(mp)^n/[n(m-1)+1]),
$$
and the prophet value $EX_n^* =P(Y^*_n=0) =g^*(0)
= 1-(mp)^n/[n(m-1)+1]$.
Note that $Y^*_n$ takes on only the two values $0$ and $n(m-1)+1$, and hence is of the same
type as the original $Y$.
\demo{Example 7.3}
 Let $X$ correspond to a Poisson ${\cal P}(\lambda)$ variable
$Y$, as in Example 4.3, with $\lambda \le 1$.
Then $k(x)=(\lambda e^{\lambda(x-1)})^n$,
so
$$
g^* (x) =(\lambda^{n-1}/n)e^{n\lambda (x-1)}+ (1-\lambda^{n-1}/n),
$$
and hence $Y_n^*$ is a mixture of a Poisson ${\cal P}(n\lambda)$ random variable with
probability $\lambda^{n-1}/n$, and the constant $0$ with probability
$(1-\lambda^{n-1}/n)$. Thus the prophet value $EX_n^*$ can be computed by
$$
P(Y^*_n=0) = 1-{\lambda^{n-1}\over n} (1-e^{-n\lambda}).
$$
\demo{Example 7.4} Let $X$ have distribution (4.4) with $p \in [0,1/2], q=1-p, b=pq, c=p$,
and $\pi=1$. It follows that $Y$ is geometric ${\cal G}(p)$,
and $g(x) = q/(1-px)$. Hence $k(x) = (pq/(1-px)^2)^n$ and we may write
$$ g^*(x) = \frac{p^{n-1}}{(2n-1)q^{n-1}} \left( \frac{q}{1-px}\right)^{2n-1}
 + \left(1 - \frac{p^{n-1}}{(2n-1)q^{n-1}}\right).$$
Hence, $Y_n^*$ is a mixture of a sum of $2n-1$ independent ${\cal G}(p)$ variables,
that is, a negative binomial, with probability $(p/q)^{n-1}/(2n-1)$, and the constant
0 with probability $1-(p/q)^{n-1}/(2n-1)$.
Thus
 the prophet value $EX_n^*$ equals
$$
P(Y^*_n=0)=q^np^{n-1}/(2n-1) +(1-(p/q))^{n-1}/(2n-1)) .
$$
\demo{Remark 7.3} In a similar way it can also be shown that in the inhomogeneous
case, when $EY_i \le 1 $ for all $i = 1, \dots, n$, the prophet variable
$X^*_n = \max (X_1, \dots, X_n)$ again corresponds to a $Y^* \in {\cal Y}$.

\par\newpage\par\flushpar

\centerline{References}

\item{[1]}  Athreya, K.B. and Ney, P.E. (1972),
{\sl Branching Processes}, Springer-Verlag, N.Y.

\item{[2]} Chow, Y., Robbins, H., and Siegmund, D.  (1971).
{\sl Great Expectations: The theory of optimal stopping}, Houghton Mifflin.

\item{[3]} D'Souza, J.C. (1995). The extinction time of the inhomogeneous branching
 process.
In: {\sl Branching Processes}, C.C. Heyde, ed. {\sl Lecture Notes in Statistics},
 {\bf 99}, Springer Verlag, 106--117.

\item{[4]} de Haan, L., (1976). Sample extremes: an elementary introduction, {\sl Statist.
Neerlandica} {\bf 30}, 161--172.

\item{[5]} Harris, T.E. (1963). {\sl The Theory of Branching Processes},
 Springer-Verlag, Berlin.

\item{[6]} Hill, T.P. and Kertz, R.P. (1981).
Ratio comparisons of supremum and stop rule expectations.
{\sl Z. Wahrscheinlichkeitstheorie verw. Gebiete} {\bf 56}, 283--285.

\item{[7]} Hill, T.P. and Kertz, R.P.  (1982).
Comparisons of stop rule and supremum expectations of i.i.d. random variables.
{\sl Ann. Probab.} {\bf 10}, 336--345.

\item{[8]} Jagers, P. (1974).
Galton-Watson processes in varying environments. {\sl J. Appl. Prob.} {\bf 11},
174--178.

\item{[9]}  Jagers, P.  (1975). {\sl Branching Processes with Biological Applications},
 John Wiley \& Sons, Ltd., London.

\item{[10]} Jirina, M. (1976).  Extinction of non-homogeneous Galton-Watson processes.
 {\sl J. Appl. Prob.} {\bf 13}, 132--137.

\item{[11]} Karlin, S. and Taylor, H.M., (1975).  {\sl A first course in stochastic processes},
 Second Ed., Academic Press, Inc., N.Y.

\item{[12]} Keiding, N. and Nielsen, J.E. (1975).
Branching processes with varying and geometric offspring
distribution. {\sl J. Appl. Prob.} {\sl 12}, 135--141.

\item{[13]} Kennedy, D.P. and Kertz, R.P. (1991).
The asymptotic behavior of the reward sequence in the optimal stopping of
i.i.d random variables. {\sl Ann. Probab.} {\bf 9}, 329--341.

\item{[14]} Leadbetter, M.R., Lindgren, G. and Rootz\'en H., (1983).
{\sl Extremes and Related Properties of Random Sequences and Processes},
Springer-Verlag, N.Y.

\item{[15]} Slack, R.S. (1968).
A branching process with mean one and possibly infinite variance,
{\sl Z. Wahrscheinlichkeitstheorie verw. Geb.} {\bf 9}, 139--145.

\bye